\documentclass[11pt,a4paper,leqno]{article}

\usepackage{epsfig}
\usepackage{amssymb}
\usepackage{amsmath}

\author{Hwajeong Kim\\
        Humboldt-Universit\"at zu Berlin\\
        Institut f\"ur Mathematik\\
        Unter den Linden 6, D-10099 Berlin\\
        Germany\\ 
        hjkim@mathematik.hu-berlin.de}
\title{A variational approach to the regularity of minimal surfaces of annulus type in Riemannian manifolds}

\usepackage{latexsym}
\usepackage{epsf,epsfig, xspace} 
\usepackage{amsmath,amsthm,amstext,amsfonts,amssymb,amsxtra,eucal} 
\newtheorem{theorem}{Theorem}[section] 
 
\newtheorem{lemma}{Lemma}[section] 
\newtheorem{remark}{Remark}[section] 
  
\newcommand{\ar}{A_{\rho}}
\newcommand{\hs}{H^{\frac{1}{2},2}} 
\newcommand{\en}{\enspace} 
\newcommand{\mhr}{\mathcal H _{\rho}}
\newcommand{\hr}{\mathcal H _{\rho}}
\newcommand{\mh}{\mathcal H}

\newcommand{\dach}{^}
\newcommand{\bqr}{\begin{eqnarray}}
\newcommand{\eqr}{\end{eqnarray}}
\newcommand{\bqrs}{\begin{eqnarray*}}
\newcommand{\eqrs}{\end{eqnarray*}}
\newcommand{\mf}{{\cal F}}
\newcommand{\mfr}{\mf_{\rho}}
\newcommand{\me}{{\cal E}}
\newcommand{\mm}{{\cal M}}

\newcommand{\mj}{{\mathbf J}}

\newcommand{\mon}{\text{mon}}

\newcommand{\dist}{\text{dist}}

\begin{document}

\maketitle

\begin{center}{\large\bf Abstract} \end{center}
Given two Jordan curves in a Riemannian manifold,   
  a minimal surface of annulus type bounded by these curves is described as 
    the harmonic extension of a critical point of some functional 
   (the Dirichlet integral) in a certain space of boundary parametrizations. 
   The $H^{2,2}$-regularity of the minimal surface of annulus type will be 
       proved by applying the critical points 
           theory and Morrey's growth condition. \medskip

{\bf Mathematics Subject Classification(2000):} 49Q05, 58E05 \smallskip

{\bf Key words} : Minimal surfaces, Abstract critical point theory, Regularity, Plateau's problem, Morrey's growth condition    

\section{Introduction}

Extending the Ljusternik-Schnirelman Theory on convex sets in Banach spaces, a general theory of critical points was developed in 1983 (\cite{s1}, see also \cite{s3} \cite{s4}), and an approach to unstable solutions and Morse theory for Plateau's problem of disc or annulus type in $\mathbb R^n$ was given. Here a minimal surface is described as the harmonic extension of a critical point of the following functional, defined on a set of boundary parametrizations:
  \[\me(x):=\frac{1}{2}\int |\mh(x)|^2d\omega,\]

where $\mh$ denotes the harmonic extension in $\mathbb R^n$. $H^{2,2}$-regularity of the above minimal surface was proved in the setting normalized by the integral condition (see \cite{s1}). In \cite{is}, further details were given and similar results were obtained for the setting normalized by the three-points condition.\medskip      

Recently, in \cite{ho}, the existence of unstable minimal surfaces of higher topological structure with one boundary in a nonpositively 
  curved Riemannian manifold was studied by applying the method introduced in \cite{s3}, and the regularity of minimal surfaces was discussed.\\

In this paper, we want to give a similar regularity result for a minimal surface of annulus type in manifolds satisfying 
  some appropriate conditions, namely, we will consider two 
  boundary curves $\Gamma_1, \Gamma_2$ in a Riemannian manifold $(N,h)$ such that one of the following conditions 
  holds.   
\begin{enumerate}
        \item[(C1)]       
                  There exists a point $p\in N$ with $\Gamma_1, \Gamma_2 \subset B(p,r)$, where  $B(p,r)$ lies within the normal range of 
                           all of its points. Here we assume $ r<\pi/(2 \sqrt{\kappa})$, where $\kappa$ is an upper bound of the 
                           sectional curvature of $(N,h)$.       
        \item[(C2)]       
                  $N$ is compact with nonpositive sectional curvature.     
  \end{enumerate} 
These conditions are related to the existence and the uniqueness of the harmonic extension for a given 
  boundary parametrization. \medskip
 
We first construct suitable spaces of functions, the sets of boundary parametrizations, where we have to distinguish 
  the cases of (C1) and (C2). 
 Then, following some idea of Struwe, we introduce a convex set 
        which, in fact, serves as a tangent space for the given boundary parametrization. 
   Moreover, we consider  the following functional: 
        \[\me(x):= \frac{1}{2}\int |d\mf(x)|_h^2d\omega, \]
where $\mf(x)$ denotes the harmonic extension of annulus type in a manifold $N$ with metric $h$.\medskip

We may then describe a minimal surface as the harmonic extension of a critical point of $\me$.\medskip

 We will always use the fact that 
  $N$ can be properly embedded into some $\mathbb R^k$ as a closed submanifold (see \cite{gr}). \medskip

Then we compute the $H^{2,2}$-regularity of our surfaces using the Morrey growth condition, see Section \ref{morreygrowth}. 
  We generalize the idea in \cite{s1} to a minimal surface of annlus type in Riemannian manifolds of the above property.

\section{Preliminaries}

\subsection{Some definitions}

Let $(M,g)$ be a manifold of dimension $2$ with boundary $\partial M$, metric $(g_{ij})$, and $(N,h)$ a connected, oriented, complete Riemannian manifold with metric $(h_{\alpha\beta})$ of dimension $n \ge 2$, embedded isometrically and properly into some $\mathbb R^k$ as a closed submanifold by $\eta$ (see \cite{gr}). Moreover, $\nabla $ and $\widetilde{\nabla}$ denote the covariant derivative in $(N,h)$ and $\mathbb R^k$, respectively.\medskip
 
   We use the summation convention for indices and a colon denotes the ordinary derivative with $i=1,2$, $\alpha = 1,\cdots,n$. Moreover, $d\omega$ and $d_0$ denote the area element in $\Omega \subset \mathbb R^2$ in $\partial \Omega$, respectively.\medskip

$\bullet$ The energy of $f\in C^2((M,g),(N,h))$ is defined by 
        \[E(f) := \frac{1}{2}\int_{M}|df|^2dM_g = \frac{1}{2} \int_{M} g^{ij}h_{\alpha\beta}\circ f 
                           f^{\alpha}_{,i}f^{\beta}_{,j}dM_g.\]

  The Euler-Lagrange equation of $E$ for $f\in C^2((M,g),(N,h))$, called the tension field along $f$, is as follows:
  \[  \tau_h(f)  :=      \langle\nabla_{\frac{\partial}{\partial z^i}}df,dz^i \rangle = g^{ij}(\nabla df)^{\alpha}_{ij}
                 =      g^{ij} ( f^{\alpha}_{,ij} -      f^{\alpha}_{,k}\Gamma^{k}_{ij}  
                                                          + f^{\beta}_{,i}f^{\gamma}_{,j}\Gamma^{\alpha}_{\beta\gamma}\circ f) 
                                                                                \frac{\partial}{\partial y^{\alpha}}\circ f. 
  \] 
        Further, $f\in  C^2((M,g),(N,h))$ is called harmonic if $\tau_h(f) = 0$.\medskip

  For $f=(f^a)_{a=1,\cdots,k}$, the second fundamental form of $\eta$ is  :
         
          \[ II\circ f (df,df) := \langle \widetilde{\nabla}_{\frac{\partial}{\partial z^i}}df
                         - \nabla_{\frac{\partial}{\partial z^i}}df, dz^i \rangle
                                                  \in T^{\bot}_{f(\cdot)}\eta(N).\]

$\bullet$ A weak Jacobi field $\mj$ with boundary $\xi$ 
   along a harmonic function $f$ 
  is a vector field along $f$ as a weak solution of 
  \[ \int_{M} \langle\nabla \mj, \nabla X \rangle + \langle tr\, R(\mj, df)df, X \rangle d\omega = 0\]
for all $X\in H^{1,2}\cap L^{\infty}(M, f^{\ast}TN)$ with $X|_{\partial M}=\xi$.\medskip

$\bullet$       For $B:=\{w\in \mathbb R^2| |w| < 1\} $,  
                 \[ H^{1,2}\cap C^0(B,N) := \{ f\in H^{1,2}\cap C^0(B,\mathbb R^k) | f(B)\subset N \}, \]
 with the norm, $ \| f \|_{1,2 ;0}:= \|\nabla f\|_{L^2} + \|f\|_{C^0}$.\medskip      
   
   Let  $\Gamma$ be a Jordan curve in $N$ that is diffeomorphic to $S^1:=\partial B$, and      
          observe that $N$ can be equipped with another metric $\tilde h$ such that $\Gamma$ is a geodesic in $(N, \tilde{h})$. 
   Note that $H^{1,2}\cap C^0\big((B,\partial B),(N,\Gamma)_{\tilde h} \big)$ and 
        $H^{1,2}\cap C^0\big((B,\partial B),(N,\Gamma)_{h} \big)$ coincide as sets. 
         Using the exponential map in $(N,\tilde{h})$, we define the following spaces.             
           \bqrs
                  \hs\cap C^0(\partial B;\Gamma) & := & \{ u\in \hs \cap C^0(\partial B, \mathbb R^k) | u(\partial B) = \Gamma \}
           \eqrs
with the norm $ \| u \|_{\frac{1}{2},2 ;0}:= \|\nabla\mh(u)\|_{L^2} + \|u\|_{C^0}$, here $\mh(u)$ is 
        the harmonic extension in $\mathbb R^k$, and  
           \bqrs
                   T_u \hs\cap C^0(\partial B;\Gamma) & := & \{ \xi \in  \hs\cap C^0(\partial B, u^{\ast}TN) | \xi(z) \in T_{u(z)}\Gamma, 
                                                                                                \en      \text{for all} \en z\in \partial B\} \\
                                                                                           & = & \hs\cap C^0(\partial B,u^{\ast}T\Gamma). 
           \eqrs

\subsection{The setting}\label{construction}

 Let $\Gamma _1, \Gamma_2$ be two Jordan curves of class $C^3$ in $N$ with diffeomorphisms 
   $\gamma^i : \partial B \rightarrow \Gamma_i, i=1,2$, and $\dist(\Gamma _1, \Gamma_2)>0$.      Moreover, for $\rho \in (0,1)$,
        \begin{eqnarray*}
                                                           A_{\rho} = \{w \in B \mid \rho < |w| < 1 \},\en 
                           C_1 = \{ w\mid|w|= 1 \},      \en C_2  = \{w \mid |w| = \rho \}. 
        \end{eqnarray*}

Let further   
         \[\mathcal X^i_{\mon}:=\{x^i\in \hs\cap C^0(\partial B ; \Gamma_i)\,|\,x^i\en\text{is weakly monotone onto}\en \Gamma_i\}.\]

 {\bf I)} We first consider the following condition for $(N,h) (\supset \Gamma_1, \Gamma_2)$:
   \begin{enumerate}
         \item[(C1)] 
                   There exists a point $p\in N$ with $\Gamma_1, \Gamma_2 \subset B(p,r)$, where  $B(p,r)$ lies within the normal range of 
                           all of its points. Here we assume $ r<\pi/(2 \sqrt{\kappa})$, where $\kappa$ is an upper bound of the 
                           sectional curvature of $(N,h)$.       
  \end{enumerate} 

 In this paper, $B(p,r)$ denotes a geodesic ball of $p\in N$ with the properties in the condition (C1).
 
\begin{remark}\label{trace}
 If $\Gamma_1, \Gamma_2\subset N$ satisfy (C1), for each $x^i\in \hs \cap C^0(\partial B;\Gamma_i)$ and $\rho\in(0,1)$ 
  there exist $g_{\rho}\in H^{1,2}\cap C^0(\overline{\ar}, B(p,r))$ and $g^i\in H^{1,2}\cap C^0(\overline{B}, B(p,r))$ with
  $g_{\rho}|_{C_1} = x^1,\, g_{\rho}|_{C_2}(\cdot) = x^2(\frac{\cdot}{\rho})$ and $g^i|_{\partial B}=x^i, \, i= 1,2$.
\end{remark}
{\bf Proof.} 
           Let $\Omega:=\exp^{-1}(B(p,r))\subset 
                          B(0,\tilde{r})_{\mathbb R^n} \subset \mathbb R^n$ for some $\tilde{r}>0$.\\ 
           For $\widetilde{x}^i:= \exp^{-1}(x^i)$, we have an Euclidean harmonic 
                 extension $h_{\rho}(\widetilde{x}^1, \widetilde{x}^2)$ of finite energy, 
                  whose image is in $B(0,\tilde{r})_{\mathbb R^n}$. 
                   The map $\exp$ is a diffeomorphism and $\Omega$ is star shaped, 
                 so there exists a retraction $\delta : B(0,\tilde{r})_{\mathbb R^n} \rightarrow \Omega$ 
                 with $\delta|_{\Omega}= Id$ in the class of $H^{1,2}$. Then the map 
                 $g_{\rho}:=\exp(\delta(h_{\rho}(\widetilde{x}^1, \widetilde{x}^2))): \ar \rightarrow \Omega$ is an 
                 $H^{1,2}\cap C^0(\overline{\ar}, B(p,r))$-extension with boundary $x^1$ and $x^2(\frac{\cdot}{\rho})$.
          We may also find an $H^{1,2}\cap C^0(\overline{B}, B(p,r))$-extension. \hfill $\Box$ \medskip

 From the results in \cite{hkw}, \cite{jk} and the above remark, we obtain a unique harmonic map of annulus and of disc type in 
  $B(p,r)\subset N$ for a given boundary mapping in the class of $\hs\cap C^0$.  
                 Now we define,  
                \bqrs 
                   M^i 
                           & := & \{x^i\in \hs\cap C^0 (\partial B;\Gamma_i)\,| x^i \en \text{is weakly monotone, orientation preserving} \}.
           \eqrs
Then $M^i$ is complete, since the $C^0$-norm preserves the monotonicity.\medskip

We now investigate another alternative condition for $(N,h)$.
 \begin{enumerate}
         \item[(C2)]   
                  $N$ is compact with nonpositive sectional curvature.     
  \end{enumerate} 
 
A compact Riemannian manifold is homogeneously regular and the condition of nonpositive sectional curvature for $N$ 
   implies $\pi_2(N) = 0$.\medskip

In order to define $M^i$, we need some preparation. First, we consider for $\rho\in(0,1)$, 
 \begin{equation*}  
         G_ {\rho} := \{f \in H^{1,2}\cap C^0(\overline{\ar},N) |\, f|_{C_i}\, \text{is continuous}      \en \text{and weakly monotone onto}
                                                                   \en \Gamma_i \}.
 \end{equation*}
We may take a continuous homotopy class, denoted by $F_ {\rho}\subset G_ {\rho}$, so that every two elements $f, g$ in $F_ {\rho}$ are continuous homotopic (not necessarily relative), denoted by $f\sim g$, 
more exactly: 
  \bqrs 
        f\sim g &\Leftrightarrow &\en \text{there exists a continuous mapping} 
                                \en H : [0,1]\times \overline{A_{\rho}} \rightarrow N \\
                 && \en \text{with} \en H(0,\cdot) = f(\cdot), H(1,\cdot) = g(\cdot). 
 \eqrs

  Now define  
                \begin{eqnarray*}
                  M^1 &:=&\{f|_{C_1}(\cdot )\in \hs\cap C^0(\partial B;\Gamma_1)\, | \,\text{orientation preserving,}\en 
                                                                  \, f \in \mathcal F_{\rho} \},\\
                  M^2 &:=& \{f|_{C_2}(\cdot \rho)\in \hs \cap C^0(\partial B;\Gamma_2)\, |\,\text{orientation preserving,}\en
                                                                  \,  f \in \mathcal F_{\rho}\}.
                \end{eqnarray*}
Then, for $x^i \in {M^i}$, there exists a unique harmonic extension to $\ar$ 
   with $x^1(\cdot)$ on $C_1$ and $x^2(\frac{\cdot}{\rho})$ on $C_2$ by \cite{le}, \cite{es}, \cite{hm}.\medskip

{\bf Definition} For $x^i\in M^i, i=1,2$, let $\mf_{\rho}(x^1,x^2)$ be the unique solution of the following Dirichlet problem:  
                 \bqr  
                   \tau_h(\mf_{\rho}(x^1,x^2)) & = & 0 \en \text{in} \en \ar \nonumber\\
                   \mf_{\rho}(x^1,x^2)(e^{i\theta})& = & x^1(e^{i\theta}) \en \text{on} \en C_1\label{dir.prob.annulus} \\
                   \mf_{\rho}(x^1,x^2)(\rho e^{i\theta}) & = & x^2(e^{i\theta}) \en \text{on}
                   \en C_2 (= \partial B_{\rho}), \nonumber 
                 \eqr 
                
  and define $ \me : \mm  \longrightarrow \mathbb R$ with 
                           \[  x \longmapsto  E(\mf(x)):=  \frac{1}{2}\int_{\ar}|d\mf_{\rho}(x^1,x^2)|_h^2d\omega.\]

{\bf II)} \,Now let $(N,h)$ and $\Gamma_i, i=1,2$, satisfy (C1) or (C2).\medskip

We will introduce a kind of tangent space of $x^i \in M^i$.\medskip

For a given oriented $y^i\in\mathcal X^i_{\mon}$, there exists a weakly monotone map $w^i\in C^0(\mathbb R, \mathbb R)$ 
  with $w^i(\theta + 2\pi) = w^i(\theta) + 2\pi$ such that 
          $y^i(\theta) = \gamma^i (\cos(w^i(\theta)), \sin(w^i(\theta)))=:\gamma^i\circ w^i(\theta)$.\\
We note that  $w^i = \tilde{w}^i + Id$ for some $\tilde{w}^i\in C^0(\partial B,\mathbb R)$. 
         Roughly speaking, $w^i$ can be viewed as a map in $C^0(\partial B, \partial B)$ 
   and then  $w^i$ is unique for given $y^i$,
        whereas $w^i\in C^0(\mathbb R, \mathbb R)$ is unique up to $2\pi l,\,l\in Z$. 
         Whether $w^i$ is in $C^0(\partial B, \partial B)$ or 
   $C^0(\mathbb R, \mathbb R)$ will be determined according to a given situation, simply denoted by  
         $y^i=\gamma^i\circ w^i$.\medskip

Denoting the Dirichlet integral by $D$ and the  $\mathbb R^k$-harmonic extension by $\mathcal H$, let 
   \begin{equation*}
         W^i_{\mathbb R^k} := \{ w^i\in C^0(\mathbb R, \mathbb R)\,|\,\text{weakly monotone},    
                           w^i(\theta + 2\pi) = w^i(\theta) + 2\pi ; 
                           D(\mathcal H (\gamma^i\circ w^i)) < \infty \}. 
        \end{equation*}
 Clearly, $W^i_{\mathbb R^k}$ is convex. For further details, we refer to \cite{s1}. \medskip
  
{\bf Definition} For $x^i\in M^i$, considering $w - w^i$ as a tangent vector along $\tilde{w}^i$, let     
        \bqrs 
         \mathcal T_{x^i} = \{ d\gamma^i((w - w^i)\frac{d}{d\theta}\circ\tilde{w}^i)
                        \,|\,w\in W^i_{\mathbb R^k}\en\text{and}\en \gamma^i\circ w\dach i =x\dach i\}.
   \eqrs

$\mathcal T_{x^i}$ is convex in $T_{x^i}\hs\cap C^0(\partial B;\Gamma_i)$, since $W^i_{\mathbb R^k}$ is 
convex.\medskip
 
Let $\widetilde{\exp}$ denote the exponential map with respect to the metric $\tilde h$. Then we note the following.  

\begin{remark}
In case of (C1), $\widetilde{\exp}_{x^i}\xi\in M^i$ for $\xi\in \mathcal T_{x^i}, i=1,2$. \\
For the case (C2), there exist $l_i>0,$ depending on $\gamma^i$ such that for any $x^i\in M^i$,
  $\widetilde{\exp}_{x^i}\xi \in M^i,\en \text{if}\en \|\xi\|_{\mathcal T_{x^i}} < l_i, i=1,2$.  
\end{remark}
{\bf Proof}      For $(C1)$ it is clear. In the case of $(C2)$, for some small $\delta>0$, 
   there exists a retraction $r$ from the $\delta$-neighborhood of 
  $N$ in $\mathbb R^k$ onto $N$, since $N$ is compact. Then, letting $\|x^i - x^i_0\|_{\frac{1}{2},2;0} < \delta$,      
                 \bqrs
                  \lefteqn { \int _{\ar}|d(r(f_\rho + \mhr(x^1 - x^1_0, 0)))|^2d\omega } \\
                        & & \le C(\|f_\rho\|_{C^0},\varepsilon,N)\big(\int_{\ar}|df_{\rho}|^2d\omega + 
                                                 \int_B|d\mathcal H(x^1-x^1_0)|^2d\omega\big)\le C(\|f_\rho\|_{1,2;0},\delta,N).
                 \eqrs
  Then we have some $l_i>0$ with the desired property, since $\widetilde{\exp}_{x^i}\xi 
                = \gamma^i(w)$ for $\xi=d\gamma^i((w - w^i)\frac{d}{d\theta}\circ\tilde{w}^i)\in\mathcal T_{x^i}$.
 \hfill $\Box$ \medskip


\begin{lemma}\label{differ}
           $\me$ is continuously partially differentiable in $x^1$ and $x^2$ with respect to variations  $ \xi^1\in\mathcal T_{x^1}$ and 
                        $ \xi^2\in\mathcal T_{x^2}$ respectively with                
                                \[ \langle \delta_{x^1}\me,\xi^1\rangle  
                                                                        =  \int_{\ar}\langle d\mfr(x^1,x^2),\nabla\mj_{\mfr}(\xi^1,0)\rangle _hd\omega.\] 
                A similar result is obtained for the second variation.  \\
        Moreover, the derivatives are continuous in $M^1\times M^2$.\medskip
          \end{lemma}
{\bf Proof} See \cite{ki}.      \hfill $\Box$ \medskip



\section{$H^{2,2}$- Regularity of minimal surfaces}

\subsection{A result}

Now we define for $x=(x^1,x^2,\rho)\in  M^1\times M^2 \times (0,1)$,
        \bqr
                  g_i(x) &:= &\sup\limits_{\begin{array}{c}\xi^i\in \mathcal T_{x^i} \\ \|\xi^i\| < l_i \end{array}}
                                                                          (-\langle \delta_{x^i}\me,\xi^i\rangle),\label{critical1} \en i=1,2.
        \eqr
Then we have the following result.
 \begin{theorem}\label{appendixthm}
        Let $x=(x^1,x^2,\rho)\in M^1\times M^2 \times (0,1)$ with $g_i(x)=0,\, i=1,2$.  Then $\mfr(x^1,x^2)$ is in the class 
         of $H^{2,2}(A_{\rho},N)$.
 \end{theorem}

\begin{remark}
In addition to the above conditions in Theorem \ref{appendixthm} let us require that 
   $g_3(x)  := \rho\cdot \partial_{\rho}\me =0$. Then, 
  $x=(x^1,x^2,\rho)$ is defined as a critical point of $\me$ such that $\mfr(x^1, x^2)$ is a minimal surface of annulus type in $N$. 
   For details we refer to \cite{ki}. 
\end{remark}

\begin{lemma}\label{quotient}
Let $\mfr := \mf_\rho(x^1,x^2) : A_{\rho} \rightarrow N \stackrel{\eta}{\hookrightarrow}\mathbb R^k$ 
          and $\mfr \in H^{1,2}(A_{\rho},\mathbb R^k)$. If      $\int_{A_\rho}|\partial_{\theta}d\mfr|^2d\omega  \le C < \infty$, then 
           $\mfr(x^1,x^2)\in  H^{2,2}(A_{\rho},N)$.
\end{lemma}
{\bf Proof} 
By Young's inequality it holds in polar coordinates with $\Delta \mfr$:=$\Delta_{\mathbb R^k} \mfr$ that 
                \bqrs
                   \lefteqn{|\nabla^{2}\mfr|^2 =  |\partial_{r}d\mfr|^2 + \frac{1}{r^2}|\partial_{\theta}d\mfr|^2}\hspace{0.2cm}\\
                                 & = & \left|\Delta \mfr - \frac{1}{r^2}\partial_{\theta\theta}\mfr - \frac{1}{r}\partial_{r}\mfr\right|^2 
                                           + \frac{1}{r^2}\left|\partial_{\theta r}\mfr\right|^2 + \frac{1}{r^4}|\partial_{\theta}\mfr|^2
                                           - 2\frac{1}{r^3}\partial_{\theta r}\mfr\partial_{\theta}\mfr
                                           +  \frac{1}{r^2}|\partial_{\theta}d\mfr|^2\\
                                & \le &  C(\varepsilon)|\Delta \mfr|^2 
                                                + (1+\varepsilon)\left|\frac{1}{r^2}\partial_{\theta\theta}\mfr + \frac{1}{r}\partial_r\mfr\right|^2
                                                 + \frac{1}{r^2}|\partial_{\theta r}\mfr|^2 + \frac{1}{r^4}|\partial_{\theta}\mfr|^2
                                           - 2\frac{1}{r^3}\partial_{\theta r}\mfr\partial_{\theta}\mfr\\
                                          & & \hspace{1.5cm} -2\varepsilon\frac{1}{r^3}\partial_{\theta r}\mfr\partial_{\theta}\mfr
                                           + \varepsilon\frac{1}{r^2}|\partial_{\theta r}\mfr|^2 
                                           + C(\varepsilon)\frac{1}{r^4}|\partial_{\theta}\mfr|^2
                                                +  \frac{1}{r^2}|\partial_{\theta}d\mfr|^2\\
                                 & \le &  C(\varepsilon)|\Delta \mfr|^2  
                                                + (2+\varepsilon)\frac{1}{r^2}|\partial_{\theta}d\mfr|^2
                                                +  C(\varepsilon)\frac{1}{r^2}\frac{1}{r^2}|\partial_{\theta}\mfr|^2\\
                                 & \le & C(\varepsilon,\eta,A_{\rho})|d\mfr|^2 + C(\varepsilon,\rho)|\partial_{\theta}d\mfr|^2,
           \eqrs
         since $\mfr$ is harmonic in $N\stackrel{\eta}{\hookrightarrow}\mathbb R^k$, i.e., $\tau_h(f)=0$.\hfill $\Box$


\subsection{The Morrey growth condition}\label{morreygrowth}

 We introduce a lemma from \cite{mo}.
 
\begin{lemma}
Let $G$ be a bounded domain in $\mathbb R^2$. Suppose $\varphi \in H^{1,2}_0(G)$, and $\psi\in L^1(G)$ satisfies 
the Morrey growth condition
  \[\int_{B_r(z_0)} |\psi| d\omega \le C_0 r^{\mu},\en  \text{for all}\en B_r(z_0).\]
Then $\psi\varphi^2 \in L^1(G)$ and for all $ B_r(z_0)$ it holds:
 \[\int_{ B_r(z_0)\cap G} |\psi\varphi^2|d\omega \le C_1C_0r^{\mu/2}\int_G |d\varphi|^2d\omega \]
for some uniform constant $C_1$.  
\end{lemma}

Let $x^i$ be as in Theorem \ref{appendixthm} with $x^i=\gamma^i\circ w^i$, and $w^i= \tilde w^i + Id$, $\tilde w^i\in H^{\frac{1}{2},2} \cap C^0(\partial B,\mathbb R)$, $i=1,2$ (recall the construction in Section \ref{construction}). Moreover, 
  for a given function $f$ on $\mathbb R$, $f_{+}(\cdot)$  and $f_{-}(\cdot)$ denote the function $f_{+}(\cdot+h)$ 
  and $f_{-}(\cdot-h)$, for $h\in \mathbb R$ respectively. \\
For $x^i \in {M^i}$ let $\mathcal H_{\rho}(x^1,x^2)$ denote the unique $\mathbb R^k$-harmonic extension 
  with boundary $x^i$ on $C_i$, $i=1,2$, and $\mathcal H(\cdot)$ the $\mathbb R^k$-harmonic extension of disc type.\\
Then we have the following growth condition.

\begin{lemma}\label{growthcondition}
   
                 For each $P_0 \in \partial A_{\rho}$ there exist $C_0, \mu, r_0 >0$ such that, for all $r\in [0,r_0]$, it holds that 
           \begin{equation}\label{sixteen}
                         \int_{A_{\rho}\cap B_r(P_0)} (|d\mfr|^2 + |d\hr(\tilde{w^1},0)|^2)d\omega \le 
                                           C_0 r^{\mu}\int_{A_{\rho}} (|d\mfr|^2 + |d\hr(\tilde{w^1},0)|^2)d\omega. 
                  \end{equation}
\end{lemma}

 \begin{remark}\label{growthremark}
   We also obtain the same result as in Lemma \ref{growthcondition} for $|d\mf_{\rho +}|^2$ (resp. $|d\mf_{\rho -}|^2$) and   
                 $|d\hr(\tilde{w^1_+},\tilde{w^2_+})|^2$ (resp. $|d\hr(\tilde{w^1_-},\tilde{w^2_-})|^2$).
 \end{remark}

 As in \cite{ho}, we observe the following.
 \begin{remark}\label{hope}
   (i)  Let $\mfr:\ar \rightarrow N$ be harmonic,       we have then  for $X\in H^{1,2}_0(A_{\rho},\mathbb R^k)$,
                \[-\int_{A_{\rho}}\langle II\circ\mfr(d\mfr\,,d\mfr),X\rangle d\omega +\int_{A_{\rho}}\langle d\mfr,dX\rangle d\omega = 0.\]
 
  (ii)     This means, for $X\in H^{1,2}(A_{\rho},\mathbb R^k)$ the above expression only depends on  the boundary of $X$. 
                Thus, for $\phi=(\phi^1,\phi^2)\in\hs\times\hs(\frac{\cdot}{\rho})$ we define    
                \bqr\label{fifty}
                   \mathbf A(\mfr)(\phi) :=      -\int_{A_{\rho}}\langle II\circ\mfr(d\mfr\,,d\mfr),X\rangle d\omega 
                         + \int_{A_{\rho}}\langle d\mfr,dX\rangle d\omega, 
                \eqr
          where $X$ is any mapping in $H^{1,2}(A_{\rho},\mathbb R^k)$ with $X|_{\partial A_{\rho}} = \phi$.\medskip

  Specially for $\phi^{i}\in H^{\frac{1}{2},2}\cap C^{0}(\partial B,(x^{i})^{\ast}T\Gamma_{i}), i=1,2$, we
                take $X:=\mj_{\mfr}(\phi^1,\phi^2)$, which is tangent to $N$ along $\mfr$, then $\langle II\circ\mfr(d\mfr\,,d\mfr),\mj_{\rho}(\phi^1,\phi^2)\rangle\equiv 0$ from the definition of the 
                 second fundamental form, so
                \bqr
                  \mathbf A(\mfr)(\phi) & = & \int_{A_{\rho}}\langle d\mfr, d\mj_{\mfr}(\phi^1,\phi^2)\rangle d\omega \label{div2} \\
                                                                & = & \int_{A_{\rho}}\langle d\mfr, d\mj_{\mfr}(\phi^1,0)\rangle d\omega  
                                                                         + \int_{A_{\rho}}\langle d\mfr, d\mj_{\mfr}(0,\phi^2)\rangle d\omega \nonumber \\
                                                                & = & \langle\partial_{x^1}\me,\phi^1\rangle 
                                                                           +  \langle\partial_{x^2}\me,\phi^2\rangle.\nonumber 
           \eqr
          Hence, for a critical point $x=(x^1,x^2,\rho)$ of $\me$,  we obtain that $\mathbf A(\mfr)(\xi)\ge 0$, 
                 for all $\xi = (\xi^1,\xi^2)\in \mathcal T_{x^1}\times \mathcal T_{x^2}$.
 \end{remark}

{\bf Proof of Lemma \ref{growthcondition}} We will show (\ref{sixteen}) in several steps.\medskip

         I) Let $P_0\in C_1$ fixed, $B_r := B_r(P_0)$, and
                   \begin{equation}
                         \tilde{w^1_0}  :=      Q^{-1} 
                          \int_{(B_{2r}\backslash B_r )\cap \partial B}\tilde{w^1} d_o,
                                                                                \en w^1_0:= \tilde{w^1_0} + Id : \mathbb R \rightarrow \mathbb R, 
                   \end{equation}
                 where $ \int_{\partial B\cap(B_{2r}\backslash B_r )}d_o         := Q $,
                   \[\tilde{\xi}_{\phi} := -\big[ \phi(|e^{i\theta}-P_0|) \big]^2(\tilde{w^1}-\tilde{w^1_0})
                                                                  \frac{\partial}{\partial\theta}\circ \bar{w^1} \in 
                                                                \hs\cap C^0 (\partial B, \bar{w^1}^{\ast}T(\partial B)),\]
   
                 where $\bar{w^1}$ means the map from $\partial B$ into itself, and $\phi\in C^{\infty}$ is a non-increasing function 
                   of $|z|$ satisfying the conditions $ 0\le \phi(z)\le 1$, $\phi\equiv 1$ if $|z| \le 2r$, 
                   $\phi\equiv 0$ if $|z|\ge 3r$, $|d \phi|\le\frac{C}{r}$, $|d^2 \phi|\le \frac{C}{r^2}$ for some $C$, fixed $r$.\medskip

                Since $(1-\phi^2)w^1+\phi^2 w^1_0\in W^1_{\mathbb R^k}$, $d\gamma^1(\tilde{\xi}_{\phi}) \in \mathcal T_{x^1}$, hence
                  \begin{equation}
                        \mathbf A(\mfr)(-d\gamma^1(\tilde{\xi}_{\phi}),0)\ge 0.
                  \end{equation} 
                Let $x^1_0 := \gamma^1(w^1_0)$, then 
                  \bqrs
                         x^1-x^1_0 & = &        d\gamma^1(w^1-w^1_0) - \int^{w^1}_{w^1_0}\int^{w^1}_{s'}d^2\gamma^1(s'')ds''ds'
                          =       d\gamma^1(w^1-w^1_0) - \alpha(w^1),
                  \eqrs
           and for small $r>0$,
                 \bqrs
                         \mathbf A(\mfr)(\phi^2(\mfr-\mfr^0)|_{C_1}, 0) & = & \mathbf A(\mfr)(\phi^2 d\gamma^1(w^1-w^1_0),0) 
                                                                                          - \mathbf A(\mfr)(\phi^2 \alpha(w^1),0)\\
                                                                           & \le & -\mathbf A(\mfr)(\phi^2 \alpha(w^1),0),
                 \eqrs
           where $\mfr^0(A_{\rho})\equiv x^1_0\in \Gamma_1$. \medskip

           On the other hand, for small $r>0, \phi^2(\mfr-\mfr^0)|_{C_2}\equiv 0$, 
            so we can take $\phi^2(\mfr-\mfr^0)$ in the 
                   definition of $\mathbf A(\mfr)$. Hence,
                   \bqrs
                        \lefteqn{ \mathbf A(\mfr)(\phi^2(\mfr-\mfr^0)|_{C_1}, 0)\hspace{3.0cm}}\\
                        \hspace{-2.0cm}  &=& \int_{A_{\rho}}\langle \phi^2d\mfr,d\mfr\rangle d\omega 
                                                 + \int_ {A_{\rho}}\langle 2\phi d\phi(\mfr-\mfr^0),d\mfr\rangle d\omega 
                                                  - \int_ {A_{\rho}}\langle \phi^2(\mfr-\mfr^0), II\circ \mfr(d\mfr,d\mfr)\rangle d\omega\\ 
                          &\le &        -\mathbf A(\mfr)(\phi^2 \alpha(w^1),0),
                  \eqrs
           and
                   \bqr
                           \lefteqn{\int_{A_{\rho}}\langle \phi^2d\mfr,d\mfr\rangle d\omega      \le   
                                  \int_ {A_{\rho}}\langle \phi^2(\mfr-\mfr^0), II\circ \mfr(d\mfr,d\mfr)\rangle d\omega}\nonumber\\ 
                                  & & \hspace{4.0cm}    - \int_ {A_{\rho}}\langle 2\phi d\phi(\mfr-\mfr^0),d\mfr\rangle d\omega 
                                                -\mathbf A(\mfr)(\phi^2 \alpha(w^1),0).\label{nine}
                   \eqr 

                For the estimate of $ -\mathbf A(\mfr)(\phi^2 \alpha(w^1),0)$, consider 
                           \[\widetilde{\star\star}:=\phi^2\int^{T^1(w^1)}_{w^1_0}\int^{T^1(w^1)}_{s'}d^2\gamma^1(s'')ds''ds'\in 
                                         H^{1,2}(A_{\rho},\mathbb R^k)\]
                with $\widetilde{\star\star}|_{C_1}= \phi^2 \alpha(w^1), \widetilde{\star\star}|_{C_2}\equiv 0$, 
                   where $w^1_0(r,\theta) = \tilde{w^1_0} + Id(r,\theta) = \tilde{w^1_0} +\theta,
                                                  \en (r,\theta)\in [\rho,1]\times \mathbb R$. \\
                By simple computation we obtain
                \bqrs
                 |\widetilde{\star\star}| & \le & C(\gamma^1,x^1) \phi^2 |\hr(\tilde{w^1},0)-\tilde{w^1_0}|^2, \\
                 |d \widetilde{\star\star}| & \le & C(\gamma^1,x^1) |\hr(\tilde{w^1},0)-\tilde{w^1_0}|^2 \phi|d\phi| 
                                           + C(\gamma^1,x^1) |d\hr(\tilde{w^1},0)||\hr(\tilde{w^1},0)-\tilde{w^1_0}|^2 \phi^2,  
                \eqrs
           and from (\ref{nine}) by Young's inequality,
         \bqrs
            \lefteqn{\int_{A_{\rho}}\langle \phi^2d\mfr,d\mfr\rangle d\omega \le  \int_  
                                 {A_{\rho}}|d\mfr|^2|\mfr-\mfr^0|\phi^2d\omega }\\
                 && \hspace{2.0cm}+ \frac{\varepsilon}{5}\int_{A_{\rho}}|d\mfr|^2 \phi^2 d\omega 
                                  + C(\varepsilon)\int_{A_{\rho}}|\mfr-\mfr^0|^2|d\phi|^2  d\omega \\
                 && \hspace{2.0cm}+  C\|\hr(\tilde{w^1},0)-\tilde{w^1_0}\|_{L^{\infty}(B_{3r})} 
                   \int_{A_{\rho}}\big( |d\mfr|^2 \phi^2 +  |\hr(\tilde{w^1},0)-\tilde{w^1_0}|^2 |d\phi|^2 \big)d\omega\\
                 && \hspace{2.0cm}+  C\|\hr(\tilde{w^1},0)-\tilde{w^1_0}\|_{L^{\infty}(B_{3r})}
                    \int_{A_{\rho}}\big( |d \hr(\tilde{w^1},0)|^2 + |d\mfr|^2 \big) \phi^2 d\omega \\
                 && \hspace{2.0cm}+  C \int_{A_{\rho}} |\hr(\tilde{w^1},0)-\tilde{w^1_0}|^2 |d\mfr|^2 \phi^2 d\omega.   
         \eqrs   
       Thus, for $r\in (0, r_0)$, sufficiently small, dependent on $\varepsilon$, $C$, 
           modulus of continuity of $\mfr-\mfr^0$ and $\hr(\tilde{w^1},0)-\tilde{w^1_0} $ we have the following estimate:
         \bqr   
           \lefteqn{\int_{A_{\rho}}\langle \phi^2d\mfr,d\mfr\rangle d\omega \le 
                 \varepsilon\int_{A_{\rho}} \big( |d\mfr|^2  +  |d\hr(\tilde{w^1},0)|^2 \big)\phi^2 d\omega}\nonumber \\
           & &  \hspace{4.0cm} +C(\varepsilon)\int_{A_{\rho}} \big( |\mfr-\mfr^0|^2 
                                         +  |\hr(\tilde{w^1},0)-\tilde{w^1_0}|^2 \big)|d\phi|^2 d\omega.\label{ten}
       \eqr

    II) We will estimate $\int_{A_{\rho}} |d\hr(\tilde{w^1},0)|^2 \phi^2 d\omega$.\medskip

      $ \bullet $ First, we obtain
            \bqrs
              \lefteqn{ D\big[ (\hr(\tilde{w^1},0)-\tilde{w^1_0})\phi \big] =
                            \int_{A_{\rho}} \big[ |d\hr(\tilde{w^1},0)|^2\phi^2 
                                     +  |(\hr(\tilde{w^1},0)-\tilde{w^1_0})|^2 |d\phi|^2}\\
              & &  \hspace{6.0cm}  + 2d\hr(\tilde{w^1},0)(\hr(\tilde{w^1},0)-\tilde{w^1_0})\phi d\phi \big] d\omega, 
            \eqrs
          and by Young's inequality
                         \bqr
                                   \lefteqn{\int_{A_{\rho}}      |d\hr(\tilde{w^1},0)|^2\phi^2 d\omega \le 
                                                   D\big[ (\hr(\tilde{w^1},0)-\tilde{w^1_0})\phi \big] }\nonumber \\
                                  && \hspace{2.0cm}      + \frac{\varepsilon}{4} \int_{A_{\rho}} |d\hr(\tilde{w^1},0)|^2\phi^2 d\omega
                                          + C(\varepsilon)\int_{A_{\rho}}\big(|\hr(\tilde{w^1},0)|^2 + |\tilde{w^1_0}|^2 \big) |d\phi|^2 d\omega.
                                   \label{thirteen}        
                        \eqr     

         $\bullet$ The estimate of $D\big[ (\hr(\tilde{w^1},0)-\tilde{w^1_0})\phi \big]$:\medskip

        On $C^1$, we have $\mfr-\mfr^0 = d\gamma^1(w^1-w^1_0) - \int^{w^1}_{w^1_0}\int^{w^1}_{s'}d^2\gamma^1(s'')ds''ds'$, and
                           $\phi|_{\partial B_{3r}(P_0)} \equiv 0$. Hence, on $\partial(A_{\rho}\cap B_{3r}(P_0))$, 
                  \bqrs
                        \lefteqn{ 
                (\hr(\tilde{w^1},0)-\tilde{w^1_0})\phi = |d\gamma^1(T^1(w^1))|^{-2} \big[ d\gamma^1(T^1(w^1))\cdot(\mfr-\mfr^0)}
                \\ 
                & & \hspace{4.0cm} 
                   + d\gamma^1(T^1(w^1))\cdot\int^{T^1(w^1)}_{w^1_0}\int^{T^1(w^1)}_{s'} 
                d^2\gamma^1(s'')ds' \big]\phi. \\
                  \eqrs 
                We denote the latter map on $\ar$ by $\Psi$.\medskip

                Moreover,  it holds that 
                   \begin{equation}\label{eleven} 
                                 \Delta\big[ (\hr(\tilde{w^1},0)-\tilde{w^1_0})\phi \big] 
                                                = 2d\hr(\tilde{w^1},0)\cdot d\phi + (\hr(\tilde{w^1},0)-\tilde{w^1_0})\Delta \phi =:f.  
                   \end{equation} 

                Note that for a solution $\varphi\in C^2(\Omega,\mathbb R)$  of  $\Delta \varphi = f$ it holds, with a 
                                  boundary data $\varphi_0$, that
                          \[D\varphi\le D\psi - \int f(\varphi-\psi),
                                 \en \text{for all}\en \psi\in \varphi_0 + H^{1,2}_0(\Omega).\]
  
                Hence, by the variation characterization of equation (\ref{eleven}), we obtain
   
                   \begin{equation}\label{twelve} 
                          D\big[ (\hr(\tilde{w^1},0)-\tilde{w^1_0})\phi \big] 
                                \le D(\Psi) - \int_{A_{\rho}\cap B_{3r}} f\big[ (\hr(\tilde{w^1},0)-\tilde{w^1_0})\phi - \Psi \big] d\omega.     
                   \end{equation}

          Let 
                        \bqrs
                                  \Psi 
                                   &:=& \frac{ d\gamma^1(T^1(w^1))\cdot(\mfr-\mfr^0)+ d\gamma^1(T^1(w^1))\cdot
                                         \int^{T^1(w^1)}_{w^1_0}\int^{T^1(w^1)}_{s'} d^2\gamma^1(s'')ds'}{|d\gamma^1(T^1(w^1))|^{2}}\phi \\
                                   & =& \frac{\Theta }{ |d\gamma^1(T^1(w^1))|^{2}}\phi,
                         \eqrs 
  
                         \[ d[d\gamma^1(T^1(w^1))\cdot(\mfr-\mfr^0)]     =      d^2\gamma^1(T^1(w^1))d(T^1(w^1))(\mfr-\mfr^0)
                                  +d\gamma^1(T^1(w^1))d\mfr=:a, \]
                         \[ d \Big( \int^{T^1(w^1)}_{w^1_0}\int^{T^1(w^1)}_{s'} d^2\gamma^1(s'')ds' d\Big) 
                                = d^2\gamma^1(T^1(w^1))d\hr(\tilde{w^1},0)(\hr(\tilde{w^1},0)-\tilde{w^1_0})
                                =:       b,\]
                         \[ d|d\gamma^1(T^1(w^1))|^{-2}  =      -2 |d\gamma^1(T^1(w^1))|^{-4}
                                \langle d^2\gamma^1(T^1(w^1)),  d^1\gamma^1(T^1(w^1)) \rangle d\hr(\tilde{w^1},0)=:c,\]
           that we have
                        \[ |d\Psi|^2 = \frac{|a+b|^2\phi^2      + \Theta^2\phi^2c^2 
                                        +  \Theta^2|d\phi|^2     + (a+b)c\phi^2\Theta   + (a+b)\phi\Theta d\phi 
                                                         +\Theta^2\phi cd\phi}{|d\gamma^1(T^1(w^1))|^{2}},\]            
                and we compute further, from the property of $\phi$, that
                         \bqrs 
                           \lefteqn{\int_{A_{\rho}} |d\Psi|^2 d\omega  \le C\int_{A_{\rho}}|d\mfr|^2\phi^2 d\omega 
                                         + C\int_{A_{\rho}}\big[ |\mfr-\mfr^0|^2 + |\hr(\tilde{w^1},0)-\tilde{w^1_0}|^2 \big]|d\phi|^2 d\omega}\\
                                  && \hspace{2.8cm}+ C\delta
                                 \int_{A_{\rho}} \big[|\hr(\tilde{w^1},0)-\tilde{w^1_0}|^2|d\phi|^2 + |d\hr(\tilde{w^1},0)|^2\phi^2 \big]d\omega,
                        \eqrs
                where $\delta=\big\||\mfr-\mfr^0|+|\hr(\tilde{w^1},0)-\tilde{w^1_0}|\big\|_{L^{\infty}(A_{\rho}\cap B_{3r})}$.\medskip

           We can also compute that 
                 \bqrs
                  \lefteqn{ -\int_{A_{\rho}\cap B_{3r}} f\big[ (\hr(\tilde{w^1},0)-\tilde{w^1_0})\phi - \Psi \big] d\omega}\\
                  && \le \int_{A_{\rho}\cap B_{3r}}\big[        2|d\hr(\tilde{w^1},0)d\phi |\hr(\tilde{w^1},0)-\tilde{w^1_0}||d \phi| 
                        + |\hr(\tilde{w^1},0) - \tilde{w^1_0}|^2|\Delta \phi|\phi  \\
                  && + C|d\hr(\tilde{w^1},0)|\phi|\mfr-\mfr^0||d\phi| 
                           +C|\hr(\tilde{w^1},0)-\tilde{w^1_0}||\mfr-\mfr^0||\Delta \phi|\phi\\ 
                  &&    +C\|\hr(\tilde{w^1},0) - \tilde{w^1_0}\|( |d\hr(\tilde{w^1},0)|\phi|\hr(\tilde{w^1},0) - 
                                \tilde{w^1_0}||d\phi| + |\hr(\tilde{w^1},0) - \tilde{w^1_0}|^2|\Delta \phi|\phi )\big] d\omega\\
                  &&    \le \int_{A_{\rho}\cap B_{3r}}\big[ C(|\mfr-\mfr^0|^2 + |\hr(\tilde{w^1},0) - \tilde{w^1_0}|)^2(|d\phi|^2+|\Delta \phi|)\\
                  && \hspace{4.0cm}      + (\frac{\varepsilon}{2} + C\|\hr(\tilde{w^1},0) - \tilde{w^1_0}\|_{L^{\infty}(A_{\rho}\cap B_{3r})})
                                                                |d\hr(\tilde{w^1},0)|^2\phi^2 \big] d\omega.  
                \eqrs
          Now the estimate of $D\big[ (\hr(\tilde{w^1},0)-\tilde{w^1_0})\phi \big]$ follows from (\ref{twelve}).\medskip

         $\bullet$      From (\ref{thirteen}) and the above estimates, we derive 
           \bqr
                 \lefteqn{\int_{A_{\rho}} |d\hr(\tilde{w^1},0)|^2 \phi^2 d\omega \le  C\int_{A_{\rho}} |d\mfr|^2\phi^2d\omega}\hspace{1.5cm}
                          \nonumber\\
                        && + C(\varepsilon)\int_{A_{\rho}}(|\mfr-\mfr^0|^2 + |\hr(\tilde{w^1},0)-\tilde{w^1_0}|^2)(|d\phi|^2+|\Delta\phi|)d\omega
                                                                        \nonumber\\
                        && +    (\frac{3\varepsilon}{4} + C\big\||\mfr-\mfr^0|+|\hr(\tilde{w^1},0)-\tilde{w^1_0}|\big\|_{L^{\infty}(A_{\rho}\cap B_{3r})})
                                          \int_{A_{\rho}}|d\hr(\tilde{w^1},0)|^2\phi^2 d\omega.\label{fourteen}
          \eqr             
 
   III) From (\ref{ten}), (\ref{fourteen}), for $r\le r_0$, where $r_0$ is dependent on $\varepsilon,\, C(x^1,\rho)$ and  the modulus of 
           continuity of $ \mfr-\mfr^0$ and $\hr(\tilde{w^1},0) -\tilde{w^1_0}$,  the definition of $\phi$\ yields
                 \bqrs
         \lefteqn{      \int_{A_{\rho}\cap B_{3r}} (|d\mfr|^2+ |d\hr(\tilde{w^1},0)|^2)d\omega 
                         \le  Cr^{-2} \int_{A_{\rho}\cap B_{3r}\backslash B_{2r}} (|\mfr-\mfr^0|^2 + |\hr(\tilde{w^1},0) - \tilde{w^1_0}|^2 )d\omega}\hspace{3.9cm} \\
                                 &\le& Cr^{-2} \int_{A_{\rho}\cap B_{3r}\backslash B_{r}}  
                                                                           (|\mfr-\mfr^0|^2 + |\hr(\tilde{w^1},0) - \tilde{w^1_0}|^2 )d\omega\\
           \text{(Poincar\'e inequality)} & \le& C \int_{A_{\rho}\cap B_{3r}\backslash B_{r}}\big( |d\mfr|^2 + |d\hr(\tilde{w^1},0)|^2 \big)d\omega\\ 
            &&\hspace{-2.0cm}   + Cr^{-2}\Big(\int_{\partial B\cap B_{2r}\backslash B_{r}}(\mfr-\mfr^0)d_o\Big)^2 
                +Cr^{-2}\Big(\int_{\partial B\cap B_{2r}\backslash B_{r}}
              \big(\hr(\tilde{w^1},0) - \tilde{w^1_0} \big)d_o\Big)^2,
         \eqrs

       where the last term is $0$ from the definition of $\tilde{w^1_0}$.\medskip

       On $\partial B$, we have 
        \[ \mfr-\mfr^0 = d\gamma^1(w^1_0)(\tilde{w^1}-\tilde{w^1}_0) + 
         \int_{\partial B\cap B_{2r}\backslash B_{r}} \int^{w^1}_{w^1_0}\int^{s'}_{w^1_0}d^2\gamma^1(s'')ds''ds',\]
       so, from the estimate in the integration and by the second inequality in Lemma \ref{poin},
         \bqrs
           \lefteqn{\int_{\partial B\cap B_{2r}\backslash B_{r}}(\mfr-\mfr^0)d_o\hspace{1.5cm} }\\
             & =&      \int_{\partial B\cap B_{2r}\backslash B_{r}}d\gamma^1(w^1_0)(\tilde{w^1}-\tilde{w^1}_0) d_o 
                + \int_{\partial B\cap B_{2r}\backslash B_{r}} \int^{w^1}_{w^1_0}\int^{s'}_{w^1_0}d^2\gamma^1(s'')ds''ds'\\ 
             & \le & C\int_{\partial B\cap(B_{2r}\backslash B_r )}|w^1-w^1_0|^2d_o \\
             & \le & Cr\int_{B\cap(B_{2r}\backslash B_r )}|d\hr(\tilde{w^1},o)|^2d\omega  
              + \frac{C}{r}\Big(\int_{\partial B\cap B_{2r}\backslash B_{r}}(\tilde{w^1}-\tilde{w^1}_0) d_o \Big)^2.
       \eqrs
       Here, the last term is again zero by the definition of  $\tilde{w^1_0}$.\medskip  
  
      Thus, 
        \bqrs
              \lefteqn{  Cr^{-2}\Big(\int_{\partial B\cap B_{2r}\backslash B_{r}}(\mfr-\mfr^0)d_o\Big)^2}\\ 
            && \le C \Big( \int_{B\cap(B_{2r}\backslash B_r )}|d\hr(\tilde{w^1},0)|^2d\omega \Big)^2 
                 \le C(x^1,\rho)\int_{B\cap(B_{2r}\backslash B_r )}|d\hr(\tilde{w^1},0)|^2d\omega,
       \eqrs    
     hence
           \[ \int_{A_{\rho}\cap B_r} \big( |d\mfr|^2+|d\hr(\tilde{w^1},0)|^2 \big)d\omega 
                 \le C\int_{A_{\rho}\cap B_{3r}\backslash B_r} \big( |d\mfr|^2+|d\hr(\tilde{w^1},0)|^2 \big)d\omega. \]

     Let $ \Upsilon(r):= \int_{ A_{\rho}\cap B_r(P_0)} \big( |d\mfr|^2+|\hr(\tilde{w^1},0)|^2 \big)d\omega$, 
       then the above inequality means that
         \[ \Upsilon(r) \le C\big(\Upsilon(3r) -     \Upsilon(r) \big),  \]
     where $C$ is independent of $r\le r_0$, for some small $r_0$.\medskip

     Then the inequality (\ref{sixteen}) follows from the Iteration-lemma. \hfill $\Box$ 



\subsection{The proof of the main theorem} 

We will give here the proof of  Theorem \ref{appendixthm}.
We begin with Poincar\'e inequality      as follows (see \cite{s1} Lemma 5.5):
\begin{lemma}\label{poin}
Let $z_0\in \partial A_{\rho}$, $B_r:=B_r(z_0)$, $G_r:= A_{\rho}\cap (B_{3r}\backslash B_r)$, 
          $K_r:= A_{\rho}\cap (B_{2r}\backslash B_r)$
           and $S_r:= \partial A_{\rho}  \cap B_{2r}\backslash B_r$. Then, for some small $r_0>0$, there exists a uniform constant $C$ 
          independent of $z_0$ such that for all $r\le r_0$ and for each $\varphi\in H^{1,2}(G_r)$:
          \bqrs 
           \int_{G_r} |\varphi|^2d\omega &\le & Cr^2\int_{G_r}|d\varphi|^2d\omega + C\left(\int_{S_r} \varphi\, d_o \right)^2, 
                 \en  \text{and} \\
                 \int_{S_r}|\varphi|^2d_o & \le& Cr \int_{K_r} |d\varphi|^2d\omega + \frac{C}{r}\left(\int_{S_r}\varphi\, d_o \right)^2,
          \eqrs
 where $d_o $ is the one-dimensional area element.
\end{lemma} 
{\bf Proof.} 
        Let $z_0, r $ be fixed. Suppose by contradiction that for a sequence ${\varphi_m}\in H^{1,2}(G_{r})$ 
         \[1\equiv\int_{G_r} |\varphi_m|^2d\omega \ge mr^2\int_{G_r}|d\varphi_m|^2d\omega+m\left(\int_{S_r} \varphi_m\, d_o \right)^2.\]
        Then $\{\varphi_m\}$ is bounded in $ H^{1,2}(G_r)$ and some subsequence, denoted again by $\{\varphi_m\}$, converges weakly to 
          some $\varphi$ in $H^{1,2}(G_r)$ but strongly in $L^2(G_r)$ by Rellich-Kondrakov. 
          From the above assumption, $d\varphi_m\rightarrow 0$ strongly.\\ 
        Thus, $\{\varphi_m\}$  converges strongly to some constant $C$ in $ H^{1,2}(G_r)$ 
          and  $\varphi_m \rightarrow C$ in $ L^2(S_r)$.\\
   On the other hand,  $\int_{S_r} \varphi_m\, d_o \rightarrow 0$, so $\varphi \equiv 0$ in $G_r$, 
         contradicting the assumption, since $\varphi_m\rightarrow \varphi$ in $L^2$.\medskip
  
   The second inequality can be proved similarly, supposing by contradiction that
          \[ 1 \equiv\int_{S_r}|\varphi_m|^2d_o \ge mr \int_{K_r} |d\varphi_m|^2d\omega + 
                                         \frac{m}{r}\left(\int_{S_r}\varphi_m\, d_o \right)^2\]
   and applying the above result for $\int_{K_r} |\varphi_m|^2d\omega$.\medskip

   By scaling, one can see that $C$ is independent of $z_0$,$r$. \hfill $\Box$\bigskip 
 
 {\bf Proof of  Theorem \ref{appendixthm}} 
           \medskip

          From Lemma \ref{quotient} and by a well known result in \cite{gt} it suffices to show that 
        \begin{equation} \label{endlich}
           \int_{A_\rho}|\Delta_{h}d\mfr|^2d\omega      \le C < \infty, 
         \end{equation}
         where $\Delta_{h}d\mfr:= \frac{d\mfr(r,\theta+h)-d\mfr(r,\theta)}{h}, h\not= 0$, and $C$ is independent of $h$. \medskip

We show (\ref{endlich}) in several steps.  The same notations as in the preceding sections will be used. \medskip

  {\bf (I)} 
%
           With $\Delta_{-h}\Delta_{h}\mfr|_{\partial B} = \Delta_{-h}\Delta_{h}\gamma^1\circ e^{iw^{1}}$ and 
                 $\Delta_{-h}\Delta_{h}\mfr|_{\partial B_{\rho}}(\cdot\rho) = \Delta_{-h}\Delta_{h}\gamma^2\circ e^{iw^{2}(\cdot)}$,
                \bqrs
                  \lefteqn{\int_{A_{\rho}}|\Delta_hd\mfr|^2d\omega      = 
                                          -\int_{A_{\rho}}\langle d\mfr,d\Delta_{-h}\Delta_{h}\mfr\rangle d\omega}\hspace{1.0cm}\\
                                 & = &-\int_{A_{\rho}}\langle II\circ \mfr(d\mfr,d\mfr),\Delta_{-h}\Delta_h\mfr\rangle d\omega
                                                                                          -\mathbf A(\mfr)(\Delta_{-h}\Delta_{h}\mfr|_{\partial A_{\rho}}).
                \eqrs                                                                                                                                                                                                                                                             
         
        Denoting $\gamma^1\circ  e^{iw^1}$ and $\gamma^2\circ e^{iw^2}$ by $\gamma^i(w^i(\theta))$, further $w^i(\cdot + h)$ and 
         $w^i(\cdot - h)$ by $w^i_{+}$ and $w^i_{-}$ respectively, we have:
                \bqrs
                   \lefteqn{\Delta_{-h}\Delta_h\gamma^i(w^i) = \Delta_{-h}\left[\frac{\gamma^i(w^i_{+}) - 
                                                                                                                  \gamma^i(w^i_{-})}{h}\right]}\\
                           & = & \Delta_{-h}\left[d\gamma^i(w^i)\left(\frac{w^i_{+}-w^i_{-}}{h}\right)+
                                                         \frac{1}{h}\int^{w^i_{+}}_{w^i}\int^{s'}_{w^i}d^2\gamma^i(s'')ds''ds'\right]\\
                           & = & d\gamma^i(w^i)(\Delta_{-h}\Delta_{h}w^i)-\frac{1}{h}\int^{w^i_{-}}_{w^i}d^2\gamma^i(s')ds'\cdot\Delta_{h}w^i_{-}
                                                   + \Delta_{-h}\left(\frac{1}{h}\int^{w^i_{+}}_{w^i}\int^{s'}_{w^i}d^2\gamma^i(s'')ds''ds'\right)\\
                           & = & d\gamma^i(w^i)(\Delta_{-h}\Delta_{h}w^i) + P^i.
                \eqrs
          Since $\gamma^i$ is smooth,  clearly 
                   $d\gamma^i(w^i)(\Delta_{-h}\Delta_{h}w^i)\in H^{\frac{1}{2},2}\cap C^{0}(\partial B,(x^{i})^{\ast}T\Gamma_{i})$.\medskip

          We write $w^i= \tilde w^i + Id$ for some $\tilde w^i\in H^{\frac{1}{2},2} \cap C^0(\partial B,\mathbb R)$ and define
                 a real valued map of $(r,\theta)\in [\rho,1]\times \mathbb R$ as follows: for $i=1$, 
           \[ T^1(w^1)(r,\theta) := H_{\rho}(\tilde w, 0)(r,\theta) + Id(r,\theta) \en \text{with} \en Id(r,\theta) = \theta, \]
          where $H_{\rho}(\tilde w, 0)$ is the harmonic extension to $A_{\rho} \approx [\rho,1]\times \mathbb R/2\pi$ with $\tilde w$ on 
                 $\partial B$ and $0$ on $\partial B_{\rho}$. Then it holds that 
           \[ T^1(w^1)(r,\theta+2\pi) = T^1(w^1)(r,\theta) + 2\pi, \en \text{for} \en (r,\theta)\in [\rho,1]\times \mathbb R,\]
                  and $e^{iT^1(w^1)}$ can be considered as a map from $\partial B$ into itself.\medskip

          Now define a map $S(P^1,0)(\cdot):A_{\rho}\rightarrow \mathbb R^k$ with the boundary $P^1$ 
           (resp. $0$) on $C_1$ (resp. $C_2$) 
                 as follows:
           \bqrs
                 S(P^1,0)(\cdot)&:=& 
                          -\frac{1}{h}\int^{T^1(w^1_{-})(\cdot)}_{T^1(w^1)(\cdot)} d^2\gamma^1(s')ds'\cdot 
                                           H_{\rho}(\Delta_h w^1_{-},0)(\cdot)\\ 
                   &&  \hspace{2cm}+\Delta_{-h}\left(\frac{1}{h}
                        \int^{T^1(w^1_{+})(\cdot)}_{T^1(w^1)(\cdot)}\int^{s'}_{T^1(w^1)(\cdot)} d^2\gamma^1(s'')ds''ds'\right).
                \eqrs
          Similarly, a map $S(0,P^2)(\cdot):A_{\rho}\rightarrow \mathbb R^k$ with the boundary $0$ (resp. $P^2$) on $C_1$ (resp. $C_2$):  
                \bqrs
                   S(0,P^2)(\cdot) &:= &-\frac{1}{h}\int^{T^2(w^2_{-})(\cdot)}_{T^2(w^2)(\cdot)} 
                                d^2\gamma^2(s')ds'\cdot H_{\rho}(0,\Delta_h w^2_{-})(\cdot)\\ 
                  &&\hspace{2cm}  +\Delta_{-h}\left(\frac{1}{h}
                   \int^{T^2(w^2_{+})(\cdot)}_{T^2(w^2)(\cdot)}\int^{s'}_{T^2(w^2)(\cdot)} d^2\gamma^2(s'')ds''ds'\right),
           \eqrs
          where $T^2(w^2_{-})(\cdot) = H_{\rho}(0,\tilde w)(\cdot) + Id(\cdot)$, and $S(0,P^2)|_{C_1}\equiv 0, S(0,P^2)|_{C_2}(\cdot \rho) = P^2(\cdot)$.\medskip

          Clearly $S(P^1,0), S(0,P^2)\in H^{1,2}(A_{\rho},\mathbb R^k)$, so letting 
                $S(P^1,P^2):= S(P^1,0)+S(0,P^2)$, 
          we have a map in $H^{1,2}(A_{\rho},\mathbb R^k)$ with boundary $(P^1,P^2)$.\medskip

          By computation, $\frac{h^2}{2}\Delta_{-h}\Delta_{h}w^i =\frac{1}{2}(w^i_{-} +w^i_{+}) - w^i$. And 
                   $\frac{1}{2}(w^i_{-} +w^i_{+})\in W^i_{\mathbb R^k}$ which is convex. Thus, by the definition of $\mathcal T_{x^i}$, 
           \[\frac{h^2}{2}d\gamma^i(w^i)(\Delta_{-h}\Delta_{h}w^i)
                                           \in \mathcal T_{x^i},\]
          and $\gamma^i(w^i)(\Delta_{-h}\Delta_{h}w^i)$ is in $\hs$, for which $A(\mfr)$ is well defined, recall Remark \ref{hope}. \medskip

          From (\ref{fifty}) and Remark \ref{hope}, since $g^1(x)=g^2(x)=0$,
           \bqrs
                   \frac{h^2}{2} A(\mfr)\left( d\gamma^1(w^1)(\Delta_{-h}\Delta_{h}w^1),0 \right)
                        =A(\mfr)\left( \frac{h^2}{2} d\gamma^1(w^1)(\Delta_{-h}\Delta_{h}w^1),0 \right)\ge 0,
           \eqrs
          so $A(\mfr)\left( d\gamma^1(w^1)(\Delta_{-h}\Delta_{h}w^1),0 \right)\ge 0$.\medskip

          Similarly,  for the second variation,
                $A(\mfr)\left( 0,d\gamma^2(w^2)(\Delta_{-h}\Delta_{h}w^2)(\frac{\cdot}{\rho})\right)\ge 0$.\medskip

          From now on we will omit the scaling term $(\frac{\cdot}{\rho})$ for the second variation.\medskip

          Moreover, from the definition of $A(\mfr)$, clearly it follows that
                \[\mathbf A(\mfr)(\phi^1+\xi^1, \phi^2+\xi^2) =  \mathbf A(\mfr)(\phi^1,\phi^2) + \mathbf A(\mfr)(\xi^1, \xi^2),\]
          if there exist $H^{1,2}$ extensions of $( \phi^1, \phi^2)$ and $(\xi^1, \xi^2)$.\medskip

          Hence, we have that
                \bqr
                   \lefteqn{\mathbf A(\mfr)\left( d\gamma^1(w^1)(\Delta_{-h}\Delta_{h}w^i),d\gamma^2(w^2)(\Delta_{-h}\Delta_{h}w^2)\right)}
                                  \hspace{1.0cm}\nonumber\\ 
                  & = & \mathbf A(\mfr)\left( d\gamma^1(w^1)(\Delta_{-h}\Delta_{h}w^1),0 \right)
                                +\mathbf A(\mfr)\left(0,d\gamma^2(w^2)(\Delta_{-h}\Delta_{h}w^2)\right)\ge 0. \label{sixty}
           \eqr
   
          Now we can compute: 
           \bqr    
                 \lefteqn{\int_{A_{\rho}}|\Delta_hd\mfr|^2d\omega 
                           =-\int_{A_{\rho}}\langle II\circ \mfr(d\mfr,d\mfr),\Delta_{-h}\Delta_h\mfr\rangle d\omega
                                                           -\mathbf A(\mfr)(\Delta_{-h}\Delta_{h}\mfr|_{\partial A_{\rho}})}
                                                        \hspace{1.0cm}\nonumber\\
                          &=& -\int_{A_{\rho}}\langle II\circ \mfr(d\mfr,d\mfr),\Delta_{-h}\Delta_h\mfr\rangle d\omega \nonumber\\
                                & & \hspace{1.0cm} -\mathbf A(\mfr)(P^1,P^2)
                                                  -\mathbf A(\mfr)\left(d\gamma^1(w^1)(\Delta_{-h}\Delta_{h}w^1),
                                                        d\gamma^2(w^2)(\Delta_{-h}\Delta_{h}w^2)\right)\nonumber\\
                 & \le & -\int_{A_{\rho}}\langle II\circ \mfr(d\mfr,d\mfr),\Delta_{-h}\Delta_h\mfr\rangle d\omega 
                                   -\mathbf A(\mfr)(P^1,P^2)\nonumber\\
                        & = & -\int_{A_{\rho}}\langle II\circ \mfr(d\mfr,d\mfr),\Delta_{-h}\Delta_h\mfr\rangle d\omega \label{first}\\
                                & & + \int_{A_{\rho}}\langle II\circ \mfr(d\mfr,d\mfr),S(P^1,0)\rangle d\omega 
                                                 + \int_{A_{\rho}}\langle II\circ \mfr(d\mfr,d\mfr),S(0,P^2)\rangle d\omega\label{second}\\
                                & & - \int_{A_{\rho}}\langle d\mfr,dS(P^1,0)\rangle d\omega 
                                           - \int_{A_{\rho}}\langle d\mfr,dS(0,P^2)\rangle d\omega.\label{third}
                \eqr
 
         {\bf (II)} For the estimstes of the above terms we need some preparation.\medskip


          First, let $s(\tau):=\tau\mf_{\rho,+}+(1-\tau)\mfr,\, 0\le \tau\le 1$, then
           \bqrs
                   \lefteqn{|\Delta_h II\circ \mfr (d\mfr,d\mfr)| 
                                   = |\frac{1}{h}\{ II\circ \mf_{\rho,+}(\mf_{\rho,+}, \mf_{\rho,+}) - II\circ \mfr(d\mfr,d\mfr)\}|}\\
                          &=& |\frac{1}{h}\{II\circ\mf_{\rho,+}(d\mf_{\rho,+},d\mf_{\rho,+})-II\circ\mfr(d\mf_{\rho,+},d\mf_{\rho,+})\\
                                & & \hspace{2.0cm} + II\circ\mfr(d\mf_{\rho,+},d\mf_{\rho,+}) - II\circ\mfr(d\mfr,d\mfr)\}|\\
                          &=& |\frac{1}{h}\{(dII(\mfr)\cdot(\mf_{\rho,+}-\mfr)+\int^{1}_{0}\int^{t}_{0}d^2II(s(\tau))
                                                  |\mf_{\rho,+}-\mfr|^2d\tau dt)(d\mf_{\rho,+},d\mf_{\rho,+}) \\
                                & & \hspace{2.0cm} + II\circ\mfr(d\mf_{\rho,+}-d\mfr,d\mf_{\rho,+}) 
                                                        + II\circ\mfr(d\mfr,d\mf_{\rho,+}-d\mfr)\}|\\
                          &=& |dII(\mfr)\cdot\Delta_h \mfr(d\mf_{\rho,+},d\mf_{\rho,+})
                                          + \frac{1}{h}\int^{1}_{0}\int^{t}_{0}d^2II(s(\tau))
                                                |\mf_{\rho,+}-\mfr|^2d\tau dt(d\mf_{\rho,+},d\mf_{\rho,+})\\
                          & &\hspace{2.0cm} + II\circ\mfr(\Delta_h d\mfr,d\mf_{\rho,+}) 
                                                  + II\circ\mfr(d\mfr,\Delta_h d\mfr)|\\
                         &\le& C(\|\mfr\|_{C^0(\ar)})[|\Delta_h \mfr||d\mf_{\rho,+}|^2 + |\Delta_h d\mfr|(|d\mf_{\rho,+}| + |d\mfr|)].
           \eqrs 
   
          Now let
                \[-\frac{1}{h}\int^{T^1(w^1_{-})}_{T^1(w^1)}d^2\gamma^1(s')ds':= \star  \quad  \text{and}  \quad
                          \frac{1}{h}\int^{T^1(w^1_{-})}_{T^1(w^1)}\int^{s'}_{T^1(w^1)}d^2\gamma^1(s'')ds''ds':=\star\star,\]
          then we have                                 
                \[ |\star| \le C(\gamma^1)|H_{\rho}(\Delta_{-h}w^1,0)|,  \quad  |\star\star|\le C(\gamma^1)|H_{\rho}(\Delta_{h}w^1,0)|,\]
          and
                \bqrs
                 |d\star| & = & \Big|-\frac{1}{h}\big[d^2\gamma^1(T^1(w^1_{-}))dT^1(w_{-}^1) -  d^2\gamma^1(T^1(w^1))dT^1(w^1)\big]\Big|\\
                                  & = & \Big|-\frac{1}{h}\big[\frac{d^2\gamma^1(T^1(w^1_{-}))- d^2\gamma^1(T^1(w^1))}
                                                                           {T^1(w^1_{-})- T^1(w^1)}\big(T^1(w^1_{-})- T^1(w^1)\big)dT^1(w_{-}^1)\\
                                  &        & \hspace{4.0cm}      
                           + d^2\gamma^1(T^1(w^1))\big(dT^1(w_{-}^1)- dT^1(w^1)\big)\big]\Big| \\
                                 & \le & C(\|\gamma^1\|_{C^3})\big(|H_{\rho}(\Delta_{-h}w^1,0)||dH_{\rho}(w^1_{-},0)| + |dH_{\rho}(\Delta_{-h}w^1,0)|\big),
                \eqrs

                \bqrs
                   |d\star\star|& = &\Big | d\big[ \frac{1}{h}\big(\int^{T^1(w^1_{+})}_{T^1(w^1)}d\gamma^1(s')ds' 
                                                                - \int^{T^1(w^1_{+})}_{T^1(w^1)}d\gamma^1(T^1(w^1))ds'\big) \big] \Big| \\  
                        & = & \Big|\frac{1}{h}\big[\frac{d\gamma^1(T^1(w^1_{+}))- d\gamma^1(T^1(w^1))}
                                   {T^1(w^1_{+})- T^1(w^1)}\big(T^1(w^1_{+})- T^1(w^1)\big)dT^1(w_{+}^1)\\
                            & & \hspace{4.0cm}  - d^2\gamma^1(T^1(w^1))dT^1(w^1)\big(dT^1(w_{+}^1)- dT^1(w^1)\big)\big]\Big| \\
                       & \le & C(\|\gamma^1\|_{C^2})|H_{\rho}(\Delta_{h}w^1,0)|\big(|dH_{\rho}(\tilde{w^1}_{+},0)| 
                                + |dH_{\rho}(\tilde{w^1},0)|\big). 
       \eqrs   

      Using the above results, we estimate (\ref{first}), (\ref{second}),(\ref{third}) for some $C\in \mathbb R$,  independent of $h$.\medskip

       First, 
         \bqrs
           (\ref{first}) & \le &  \int_{A_{\rho}}\big|\langle \Delta_h II\circ\mfr(d\mfr,d\mfr),\Delta_h \mfr\rangle \big| d\omega \\ 
                 & \le & C\int_{A_{\rho}}\big( |\Delta_h\mfr|^2|d\mf_{\rho,+}|^2 + 
                               |\Delta_h d\mfr|(|d\mf_{\rho,+}|+|d\mfr|)|\Delta_h \mfr|\big)d\omega \\
                 & \le &  C\int_{A_{\rho}}|d\mf_{\rho,+}|^2|\Delta_h \mfr|^2d\omega  +\varepsilon\int_{A_{\rho}}| \Delta_h d\mfr|^2d\omega + 
                             C(\varepsilon)\int_{A_{\rho}}(|d\mf_{\rho,+}|^2+|d\mfr|^2)|\Delta_h\mfr|^2d\omega.
        \eqrs 

      For the estimate of (\ref{second}),
         \bqrs 
           \lefteqn{\int_{A_{\rho}}\langle II\circ \mfr(d\mfr,d\mfr),S(P^1,0)\rangle d\omega}\\
                & \le & \int_{A_{\rho}}\{|\langle II\circ \mfr(d\mfr,d\mfr),
                            (\star) H_{\rho}(\Delta_{h}w^1_{-},0)\rangle| 
                  + |\langle\Delta_{h} II\circ\mfr(d\mfr,d\mfr), (\star\star) \rangle |\}d\omega\\
               & \le & C\int_{A_{\rho}}|d\mfr|^2|\hr(\Delta_{-h}w^1,0)|^2d\omega \\
                   & & + C\int_{A_{\rho}}\{|\Delta_{h}\mfr||d\mf_{\rho,+}|^2|\hr(\Delta_{h}w^1,0)|
                                   +|\Delta_h d\mfr|(|d\mf_{\rho,+}|+|d\mfr|)|\hr(\Delta_{h}w^1,0)|\}d\omega\\
               & \le & C  \int_{A_{\rho}}|d\mfr|^2|\hr(\Delta_{-h}w^1,0)|^2d\omega 
                 + C \int_{A_{\rho}}|d\mf_{\rho,+}|^2(|\Delta_h \mfr|^2+ |\hr(\Delta_{h}w^1,0)|^2)d\omega \\
              & & +  \varepsilon\int_{A_{\rho}}|\Delta_h d\mfr|^2d\omega + 
                C(\varepsilon)\int_{A_{\rho}}(|d\mf_{\rho,+}|^2|\hr(\Delta_{h}w^1,0)|^2+|d\mfr|^2|\hr(\Delta_{h}w^1,0)|^2)d\omega,
        \eqrs
    note that $\Delta_{h}w^1_{-} = \Delta_{-h}w^1$, and we obtain a similar estimate for the second term of (\ref{second}).\\
    Thus, we have that
       \bqrs
          (\ref{second}) & \le &  \varepsilon C\int_{A_{\rho}}|\Delta_h d\mfr|^2d\omega   
            + C(\varepsilon)\int_{A_{\rho}}\big(|d\mfr|^2 + |d\mf_{\rho,+}|^2\big)\cdot\\ 
            && \big(|\Delta_h\mfr|^2 + |\hr(\Delta_{-h}w^1,0)|^2+ |\hr(0,\Delta_{-h}w^2)|^2    +\hr(\Delta_{h}w^1,0)|^2+|\hr(0,\Delta_{h}w^2)|^2 \big)d\omega.
      \eqrs

    For the estimate of (\ref{third}), 
       \bqrs
        \lefteqn{-\int_{A_{\rho}}\langle d\mfr,dS(P^1,0)\rangle d\omega \le
          \int_{A_{\rho}}\big|\langle d\mfr, d(\star)\hr(\Delta_{-h}w^1,0)\rangle\big| d\omega} \\
          & & \hspace{2.5cm} + \int_{A_{\rho}}\big|\langle d\mfr, (\star)d\hr(\Delta_{-h}w^1,0)\rangle\big| d\omega  
           + \int_{A_{\rho}}\big|\langle \Delta_h d\mfr, d(\star\star)\rangle\big| d\omega \\
          & \le & \varepsilon C\int_{A_{\rho}}|\Delta_h d\mfr|^2d\omega + 
                                  \varepsilon C\int_{A_{\rho}}|d\hr(\Delta_h w^1,0)|^2d\omega \\
           &  &\hspace{0.5cm} + C(\varepsilon)\int_{A_{\rho}} \big(|d\mfr|^2 + |d\hr(\tilde{w}^1_{-},0)|^2+ 
                               |d\hr(\tilde{w}^1_{+},0)|^2+|d\hr(\tilde{w}^1,0)|^2\big)\cdot \\
                              & &  \hspace{5.0cm}\big(|\hr(\Delta_{-h}w^1,0)|^2+ |\hr(\Delta_{h}w^1,0)|^2\big)d\omega.
       \eqrs
    We obtain a similar estimate for the second term of (\ref{third}):  
       \bqrs 
            (\ref{third}) & \le &  \varepsilon C\int_{A_{\rho}}|\Delta_h d\mfr|^2d\omega
                         + \varepsilon C\int_{A_{\rho}}|d\hr(\Delta_h w^1,0)|^2d\omega\\
           & & + C(\varepsilon)\int_{A_{\rho}}\big(|d\mfr|^2 +
                                   |d\hr(\tilde{w}^1_{-},0)|^2+ |d\hr(\tilde{w}^1_{+},0)|^2+|d\hr(\tilde{w}^1,0)|^2\\
          & &\hspace{4.5cm} + |d\hr(0,\tilde{w}^2_{-})|^2+ |d\hr(0,\tilde{w}^2_{+})|^2+|d\hr(0,\tilde{w}^2)|^2   \big)\cdot \\
          & &  \hspace{0.8cm}\big(|\hr(\Delta_{-h}w^1,0)|^2+ |\hr(\Delta_{h}w^1,0)|^2
                                                  +|\hr(0,\Delta_{-h}w^2)|^2+ |\hr(0,\Delta_{h}w^2)|^2\big)d\omega.
      \eqrs

    Now, gathering all the above results we obtain :
       \bqr
         \int_{A_{\rho}}|\Delta_h d\mfr|^2 d\omega 
          =\varepsilon C\int_{A_{\rho}}|\Delta_h d\mfr|^2d\omega
                         + \varepsilon C\int_{A_{\rho}}|d\hr(\Delta_h w^1,0)|^2d\omega +C(\varepsilon)\Xi \,, \label{fourth}
       \eqr
    where
       \bqrs
          \Xi &:= &\int_{A_{\rho}}\big(|d\mfr|^2 + |d\mfr|^2 
	        + |d\hr(\tilde{w}^1_{-},0)|^2+ |d\hr(\tilde{w}^1_{+},0)|^2+|d\hr(\tilde{w}^1,0)|^2\nonumber\\
          & &\hspace{5.3cm} 
	       +|d\hr(0,\tilde{w}^2_{-})|^2+ 
                    |d\hr(0,\tilde{w}^2_{+})|^2+|d\hr(0,\tilde{w}^2)|^2   \big)\cdot\nonumber \\
          & & 
	      \big(|\Delta_h\mfr|^2+|\hr(\Delta_{-h}w^1,0)|^2+ |\hr(\Delta_{h}w^1,0)|^2
		     +|\hr(0,\Delta_{-h}w^2)|^2+ |\hr(0,\Delta_{h}w^2)|^2\big)d\omega\nonumber\\
        \eqrs

{\bf (III)} On  $\partial B$, it holds that 
         $ \Delta_h(\gamma^i\circ w^i) = d\gamma^i(w^i)\Delta_hw^i + 
                  \frac{1}{h}\int^{w^i_{+}}_{w^i}\int^{s'}_{w^i}d^2\gamma^i(s'')ds''ds' $, 
       so  
          \begin{equation}\label{fifth} 
               \Delta_h w^i = |d\gamma^i(w^i)|^{-2}\big[ d\gamma^i(w^i)\cdot \Delta_h\mfr 
                - d\gamma^i(w^i)\cdot\frac{1}{h}\int^{w^i_{+}}_{w^i}\int^{s'}_{w^i}d^2\gamma^i(s'')ds''ds'\big]. 
          \end{equation}

       Using $T^i(w^i)$ at the right hand side of (\ref{fifth}), we obtain an $H^{1,2}(A_{\rho},\mathbb R^k)$- extension 
           with boundary $\Delta_h w^i$ on $C^1$ and $0$ on $C_2$, and by the D-minimality of the harmonic extension 
           between the maps with the same boundary, we have
        \bqr
           \lefteqn{\int_{A{\rho}}|d\hr(\Delta_h w^1,0)|^2 d\omega}\nonumber\hspace{0.7cm}\\
            & \le & C\int_{A_{\rho}} \big[ |d\hr(w^1,0)| \big( |\Delta _h\mfr|+|\star\star| \big) + 
                           |d\Delta_h \mfr| + |d\star\star| \big]^2 d\omega \nonumber\\
             & \le & C\int_{A_{\rho}} \big\{ |d\hr(w^1,0)|^2 |\Delta_h\mfr|^2 + |d\hr(\Delta_h w^1,0)|^2|\hr(\Delta_h w^1,0)|^2
                                                                                              + |d\Delta_h \mfr|^2\nonumber \\
            & & \hspace{1.0cm} + |\hr(\Delta_h w^1,0)|^2(|d\hr(\tilde{w^1}_+ ,0)| + |d\hr(\tilde{w^1} ,0)|)^2\nonumber \\
            & & \hspace{1.0cm} + |d\hr(\tilde{w^1} ,0)|^2 |\Delta_h w^1,0)| 
                                                      + |d\hr(\tilde{w^1} ,0)||\Delta_h \mfr||d\Delta_h \mfr|\nonumber \\
             & & \hspace{1.0cm} + |d\hr(\tilde{w^1} ,0)||\hr(\Delta_h w^1,0)|(|d\hr(\tilde{w^1}_+ ,0)|
                                                                   + |d\hr(\tilde{w^1} ,0)|)|\Delta_h \mfr|\nonumber\\
              & & \hspace{1.0cm} + |d\hr(\tilde{w^1} ,0)||\hr(\Delta_h w^1,0)||d\Delta_h \mfr|\nonumber\\
             & & \hspace{1.0cm} + |d\hr(\tilde{w^1} ,0)||\hr(\Delta_h w^1,0)||\hr(\Delta_h w^1,0)|
                                                         (|d\hr(\tilde{w^1}_+ ,0)|+ |d\hr(\tilde{w^1} ,0)|)\nonumber\\
             & & \hspace{1.0cm} + |d\Delta_h \mfr||\hr(\Delta_h w^1,0)|
                                                  (|d\hr(\tilde{w^1}_+ ,0)|+ |d\hr(\tilde{w^1} ,0)|) \big\}d\omega\nonumber\\
           & \le &  C\int_{A_{\rho}}|d\Delta_h \mfr|^2 d\omega  + C\Xi\, .\label{sixth}
        \eqr
    Similarly, we obtain an estimate 

          \begin{equation}\label{seventh}
                  \int_{A{\rho}}|d\hr(0,\Delta_h w^2)|^2 d\omega \le C\int_{A_{\rho}}|d\Delta_h \mfr|^2 d\omega  + C\Xi. 
          \end{equation}

          Using the estimate (\ref{fourth}) for $\int_{A_{\rho}}|d\Delta_h \mfr|^2 d\omega$ and from (\ref{sixth}), (\ref{seventh}), 
                \bqrs
                   \lefteqn{\int_{A_{\rho}}|d\Delta_h \mfr|^2 d\omega + \int_{A_{\rho}}|d\hr(\Delta_h w^1,0)|^2 d\omega
                                                                                                   + \int_{A_{\rho}}|d\hr(0,\Delta_h w^2)|^2 d\omega}\hspace{1.0cm}\\
                  & \le & \varepsilon C \int_{A_{\rho}}|d\Delta_h \mfr|^2 d\omega + \varepsilon C\int_{A_{\rho}}|d\hr(\Delta_h w^1,0)|^2 d\omega
                                                                                                   + \varepsilon C\int_{A_{\rho}}|d\hr(0,\Delta_h w^2)|^2 d\omega
                   + C(\varepsilon) \Xi\,.
          \eqrs
         Since    $ \frac{1}{2}(a^2+b^2) \le (a+b)^2 \le \frac{3}{2}(a^2+b^2), a,b\in \mathbb R \en \text{and}\en 
                          H_{\rho}(f,g) = H_{\rho}(f,0)+H_{\rho}(0,g)$, for some small $\varepsilon>0$ 
                          in the above estimate we finally obtain the following inequality:
         \bqr
                \lefteqn{\int_{A_{\rho}}|\Delta_h d\mfr|^2 d\omega + \int_{A_{\rho}}|d\hr(\Delta_h w^1,\Delta_h w^2)|^2 d\omega}
                                                                                                                                                          \nonumber\hspace{0.3cm} \\
                          & \le & C(\varepsilon) \int_{A_{\rho}} \big( |d\mfr|^2 + |d\mf_{\rho +}|^2    + |d\mf_{\rho -}|^2 \nonumber\\
                        &&\hspace{2.0cm}+ |d\hr(\tilde{w^1},\tilde{w^2}|)^2 + |d\hr(\tilde{w^1_+},\tilde{w^2_+})|^2
                                                                                                         + |d\hr(\tilde{w^1_{-}},\tilde{w^2_-}|^2 \big)\cdot \nonumber\\
                           & &\label{eight}      \hspace{2.0cm} \big( |\Delta_h \mfr|^2 + |H(\Delta_{-h}w^1,\Delta_{-h}w^2)|^2 + 
                                                                                                                           |H(\Delta_{h}w^1,\Delta_{h}w^2)|^2 \big) d\omega .
          \eqr

{\bf (IV)} Now extend $\mfr$ to $\mathbb R^2\backslash B_{\rho^2}$ by conformal reflection as follows
                \bqrs
                        \mfr(z) & = & \mfr\big( \frac{z}{|z|^2} \big) ,\en \text{if} \en 1\le |z| \\
                        \mfr(z) & = & \mfr\big( \frac{z}{|z|^2}\rho^2 \big) ,\en \text{if} \en \rho^2 \le |z|\le \rho.
                \eqrs 
         Choose $r\in \big( 0, \min\{\frac{\rho-\rho^2}{2}, r_0\} \big)$, and $ \varphi \in C^{\infty}_0\big(B_{2r}(0) \big)$
                 with $\varphi \equiv 1$ on $B_{r}(0)$.\medskip

         We may cover $A_{\rho}$ with balls of radius $r$ in such a way that at most $k$ balls of the covering intersect at any point 
                $p\in A_{\rho}$, for any $r$ as above ($\mathbb R^2$ is metrizable).  Let $B^i$ denote the balls of the covering 
                with centers $p_i$ and  $\varphi_i(p):= \varphi(p-p_i)$. \\
         Then, from (\ref{eight})
                 \bqrs                    
                        \lefteqn{\int_{A_{\rho}}|\Delta_h d\mfr|^2 d\omega + \int_{A_{\rho}}|d\hr(\Delta_h w^1,\Delta_h w^2)|^2 d\omega}\\
                        & \le & C\,\Sigma_{i} \int_{\mathbb R^2\backslash A_{\rho^2}}
                                           \big( |\Delta_h \mfr|^2 + |H(\Delta_{-h}w^1,\Delta_{-h}w^2)|^2 + 
                                                                                                                           |H(\Delta_{h}w^1,\Delta_{h}w^2)|^2 \big)\varphi_i^2 \cdot \\
                        & &        \underbrace{\big( |d\mfr|^2 + |d\mf_{\rho +}|^2      + |d\mf_{\rho -}|^2 
                                + |d\hr(\tilde{w^1},\tilde{w^2})|^2 + |d\hr(\tilde{w^1_+},\tilde{w^2_+})|^2
                                                                                                         + |d\hr(\tilde{w^1_{-}},\tilde{w^2_-})|^2 \big)}_{=: \chi}d\omega.
           \eqrs
        According to  Lemma \ref{growthcondition} and Remark \ref{growthremark}, $\chi$ satisfies the Morrey growth condition, so apply the Morrey Lemma   with $\chi$ and $ (\Delta_h \mfr)\varphi_i $ resp.$H(\Delta_{-h}w^1,\Delta_{-h}w^2)\varphi_i$ resp. $H(\Delta_{h}w^1,\Delta_{h}w^2)\varphi_i$. 
   Then we obtain
         \bqrs
                 \lefteqn{ \int_{B_{2r}(p_i)} \chi \big( |\Delta_h \mfr|^2 + |H(\Delta_{-h}w^1,\Delta_{-h}w^2)|^2 + 
                                                                                                                           |H(\Delta_{h}w^1,\Delta_{h}w^2)|^2 \big) \varphi_i^2 d\omega}\\ 
                &&      \le Cr^{\frac{\mu}{2}}\int_{B_2\backslash B_{\rho^2}}\chi d\omega 
                        \int_{B_{2r}(P_i)}\big(|d \Delta_h \mfr|^2 + |d H(\Delta_{-h}w^1,\Delta_{-h}w^2)|^2 
                                                 + |dH(\Delta_{h}w^1,\Delta_{h}w^2)|^2 \big)d\omega \\
           && + Cr^{\frac{\mu}{2}}\int_{B_2\backslash B_{\rho^2}}\chi d\omega 
                  \int_{B_{2r}(P_i)}\big(|\Delta_h \mfr|^2 + |H(\Delta_{-h}w^1,\Delta_{-h}w^2)|^2 
                                                 + |H(\Delta_{h}w^1,\Delta_{h}w^2)|^2 \big)d\omega.
         \eqrs   

        Summing over $i$ yields a constant $C$, independent of $r$, such that 
                 \bqrs
                        \lefteqn{\int_{A_{\rho}}|\Delta_h d\mfr|^2 d\omega + \int_{A_{\rho}}|d\hr(\Delta_h w^1,\Delta_h w^2)|^2 d\omega}\\ 
                         && \le Cr^{\frac{\mu}{2}}\int_{B_2\backslash B_{\rho^2}}\big(|d \Delta_h \mfr|^2 + |d H(\Delta_{-h}w^1,\Delta_{-h}w^2)|^2 
                                                 + |dH(\Delta_{h}w^1,\Delta_{h}w^2)|^2 \big)d\omega \\
                         && +    Cr^{\frac{\mu}{2}}\int_{B_2\backslash B_{\rho^2}}\big(|\Delta_h \mfr|^2 + |H(\Delta_{-h}w^1,\Delta_{-h}w^2)|^2 
                                                 + |H(\Delta_{h}w^1,\Delta_{h}w^2)|^2 \big)d\omega.
                \eqrs 
        
   Since $d\mfr,\,dH(w^1,w^2)\in L^2$, choosing small $r>0$, we obtain $C\in \mathbb R$, independent of $|h|\le h_0$ with        
        \[ \int_{A_{\rho}}|\Delta_h d\mfr|^2 d\omega  \le C. \] 
                                                                                                                                                                                           \hfill $\Box$

\end{document}